\documentclass[12pt,leqno]{article}
\usepackage[latin1]{inputenc}
\usepackage[english]{babel}
\usepackage{amsopn,amsthm,amssymb,amsmath,epsfig,graphics,enumerate,array,calc}

\theoremstyle{plain}
\newtheorem{theo}{Theorem}[section]
\newtheorem{prop}[theo]{Proposition}
\newtheorem{lem}[theo]{Lemma}
\newtheorem{cor}[theo]{Corollary}
\newtheorem{theodef}[theo]{Theorem-Definition}

\theoremstyle{definition}
\newtheorem{Def}[theo]{Definition}

\theoremstyle{remark}
\newtheorem{ex}[theo]{Example}
\newtheorem{remark}[theo]{Remark}

\newcommand{\pil}{\rightarrow}

\newcommand{\nn}{\mathbb{N}}
\newcommand{\zz}{\mathbb{Z}}
\newcommand{\pp}{\mathbb{P}}

\newcommand{\rr}{\mathbb{R}}

\newcommand{\la}{\lambda}

\newcommand{\Ga}{\Gamma}
\newcommand{\e}{\varepsilon}

\newcommand{\nsp}{\negmedspace}
\newcommand{\oo}{\mathcal{O}}

\newcommand{\skrifta}{\mathcal{A}}

\newcommand{\av}{\boldsymbol{a}}

\newcommand{\feit}[1]{\boldsymbol{#1}}
\newcommand{\sub}{\subseteq}
\newcommand{\wh}[1]{\widehat{#1}}
\newcommand{\wb}[1]{\,\overline{\!#1\!}\,}

\DeclareMathOperator{\Subdiv}{Subdiv}
\DeclareMathOperator{\Area}{Area}
\DeclareMathOperator{\divv}{div}
\DeclareMathOperator{\Div}{Div}

\DeclareMathOperator{\Jac}{Jac}

\title{The group law on a tropical elliptic curve}
\author{Magnus Dehli Vigeland \\
Department of Mathematics, University of Oslo, Norway \\
{\it Email}\,: {\tt magnusv@math.uio.no}
}
\date{}
\begin{document}\maketitle
\section{Introduction}

Tropical elliptic curves have a natural group structure in analogy to classical elliptical curves. We define the {\em Jacobian} as an abelian group associated to a tropical curve. Unlike in the classical case, the Jacobian of a tropical elliptic curve $C$ is not equal as a set to the curve itself, but to a smaller part of it, namely the complement $\wb{C}$ of the so-called {\em tentacles} of the curve. A distance function $d_C(P,Q)\colon \wb{C}\times \wb{C}\pil \rr$ plays a crucial role in the main results:

\begin{theo} Let $C$ be a tropical elliptic curve, and let $\wb{C}$ be the complement of the tentacles. Let $\oo$ be a point on $\wb{C}$. 
  \begin{enumerate}
\item We have a bijection of sets $\wb{C}\longleftrightarrow \Jac(C)$, given by $P\longleftrightarrow P-\oo$. 
\item The induced group law on $\wb{C}$ satisfies the relation $$d_C(\oo,P+Q)=d_C(\oo,P)+d_C(\oo,Q).$$
\item As a group, $\wb{C}$ is isomorphic to $S^1$.

  \end{enumerate}
\end{theo}

\section{Definitions}
Let $\rr_{tr}:=(\rr,\oplus,\odot)$ be the {\em tropical semiring}, where the binary operations are defined by
\begin{equation}\label{trdef}
  a\oplus b := \max\{a,b\}\quad\text{and}\quad a\odot b :=a+b.
\end{equation}
The multiplicative identity element of $\rr_{tr}$ is $0$, while there is no additive identity (unless we include $-\infty$ as an element of $\rr_{tr}$). 

We extend $\oplus$ and $\odot$ to $\rr^n$ by using \eqref{trdef} on each coordinate: 
\begin{equation*}
  \begin{split}
(a_1,\dotsc,a_n)\oplus (b_1,\dotsc,b_n)&:= (\max\{a_1,b_1\},\dotsc,\max\{a_n,b_n\}),\quad \text{and}\\
\la\odot(a_1,\dotsc,a_n)&:=(\la+a_1,\dotsc,\la+a_n),\quad\text{ for $\la\in\rr$.}
  \end{split}
\end{equation*}
Furthermore, let the {\em tropical projective $n$-space} be defined by
$\pp^n_{tr}:=\rr^{n+1}\big{/}\sim\;,$ 
where $\feit{x}\sim \feit{y} \Longleftrightarrow \feit{x}=\la \odot \feit{y}$ for some $\la\in\rr$. 

\begin{remark}\label{ekstra}
Note that we don't get any additional points by going from $\rr^n$ to $\pp_{tr}^n$.
For example, every equivalence class in $\pp_{tr}^n$ has a representative in $\rr^{n+1}$ with 0 as the last coordinate. Still, we will often work projectively, since this gives a symmetrization of the coordinates and makes the presentation more canonical.
\end{remark}

Let $\skrifta\sub\zz^n$ be a finite set of vectors $\av=(a_1,\dotsc,a_n)$. A tropical (Laurant) polynomial in indeterminates $x_1,\dotsc,x_n$, with support $\skrifta$, is an expression of the form
\begin{equation*}
f=\bigoplus_{\av\in\skrifta} \la_{\av}\odot x_1^{a_1}\odot\dotsb\odot x_n^{a_n}=\max_{\av\in\skrifta}\{\dotsc,\la_{\av}+ \sum_{i=1}^n a_i x_i,\dotsc \},
\end{equation*}
where each $\la_{\av}\in\rr$. Notice that as a function $\rr^n\pil \rr$, $f$ is convex and piecewise-linear. 
The tropical polynomial $f$ is called {\em homogeneous of degree $d$} if $a_1+\dotsb+a_n=d$ for all $\av\in\skrifta$. Furthermore, the convex hull of $\skrifta$ is called the {\em Newton polytope of $f$} and is denoted by $\Delta$.

\section{Tropical curves}
\subsection{Basic properties}
To define tropical plane curves, let us first define tropical varieties in general: 

\begin{Def}
Let $f$ be a tropical polynomial in $n$ indeterminates. The {\em tropical variety} $V(f)$ defined by $f$ is the set of points in $\rr_{tr}^n$ where the associated function $f\colon \rr^n\pil\rr$ is not linear. If $f$ is homogeneous, we can regard $V(f)$ as a {\em tropical projective variety} in $\pp_{tr}^{n-1}$. 
\end{Def}

\begin{remark}\label{monom}
Note that if $f$ consists of a single monomial, $V(f)$ is an empty set.
\end{remark}

Different tropical polynomials can have the same associated tropical variety. In particular, it is easy to see that if the supports of $f$ and $g$ differ only by a translation, then $V(g)=V(f)$.

\begin{Def}
By a {\em tropical curve} in $\pp_{tr}^2$, we mean a tropical projective variety of the form $V(f)$, where $f$ is a homogeneous tropical polynomial in 3 indeterminates.
\end{Def}

We will next recall some basic properties of tropical curves. For proofs and more details,
see \cite[Section 3]{RGST}, or \cite[Sections 1-3]{Mikh} for a more exhaustive approach. 

Given a tropical polynomial $f$, we can associate a lattice
subdivision of the 
Newton polygon $\Delta$ of $f$ in the following way: Let $\wh{\Delta}$
be the convex hull of the set $\{a,b,c,\la_{abc}\}\sub\zz^3\times
\rr$, where $(a,b,c)$ runs through $\skrifta$. Then define
$\Subdiv_f$ to be the image under the projection to $\zz^3$ of the corner
edges on the upper part of $\wh{\Delta}$. 

The subdivision $\Subdiv_f$ is in a natural way dual to the
tropical variety $V(f)$. In particular, each edge of $V(f)$ corresponds to
an edge of $\Subdiv_f$, and corresponding edges are perpendicular to
each other. The unbounded rays in $V(f)$ correspond to the edges of
$\partial\Delta$. (Cf. \cite[Proposition 3.5]{RGST} and \cite[Proposition 3.11]{Mikh}.)

Let $E$ be an edge of a tropical curve $C=V(f)$, and let $\Delta'$ be
the corresponding edge in $\Subdiv_f$. We define the {\em weight} of $E$ to
be the lattice length of $\Delta'$, i.e. $ 1+\text{the number of interior lattice points of
$\Delta'$}$.
 
\begin{lem}
For any node $V$ in a tropical curve $C\sub\pp_{tr}^2$, the following {\em
  balancing condition} holds: Let $E_1,\dots,E_n$ be the edges
adjacent to $V$. For each $i=1,\dotsc,n$ let $m_i$ be the weight of
$E_i$, and $\feit{v}_i$ the primitive integer vector starting at $V$ and pointing in the
direction of $E_i$. Then
\begin{equation}\label{balance}
  m_1\feit{v}_1+\dotsb+m_n\feit{v}_n=\feit{0},
\end{equation}
where $\feit{0}=(0,0,0)\in \pp_{tr}^2$.
\end{lem}

The balancing condition characterizes the tropical curves: Assume $C$ is a collection of rays and line segments in $\pp_{tr}^2$, all with rational slopes, and all assigned some positive integral weight. Then $C=V(f)$ for some homogeneous tropical polynomial $f$, if and only if \eqref{balance} is satisfied at every vertex of $C$.

Next we define the degree of a
tropical curve. For each $d\in\nn_0$, let $\Ga_d$ be the triangle in $\zz^3$ with vertices $(d,0,0),(0,d,0),(0,0,d)$. (When $d=0$ we get the degenerated triangle $\Ga_0=\{0,0,0\}$.)

\begin{Def}\label{degdef}
Let $C=V(f)$ be a tropical curve in $\pp_{tr}^2$, and let $\Delta$ be
the Newton polygon of $f$. If $\Delta$ fits inside $\Ga_d$, but not
inside $\Ga_{d-1}$, then $C$ has {\em degree} $d$. If $\Delta=\Ga_d$,
we say that $C$ has degree $d$ {\em with full support}.
\end{Def}

\begin{remark}
There seems to be no clear consensus in the literature on how to define
the degree of a tropical curve. Definition \ref{degdef} differs slightly from the ones in \cite{RGST} and \cite{Mikh}, but serves the purpose of this paper better. In particular, as we will see in the next section, Definition \ref{degdef} gives room for an extended version of the tropical Bezout's theorem compared to that in \cite{RGST}. 
\end{remark}

\begin{ex}
A {\em tropical line} is a tropical curve of degree 1. For instance, if $f$ is the degree 1 polynomial $f=a  x\oplus b y\oplus c z$ (with Newton polygon $\Delta=\Ga_1$), then $L=V(f)$ is the tropical line in with ''center'' $(-a,-b,-c)\in\pp_{tr}^2$. Note that the same line could be described as $V(a x^2\oplus b xy\oplus c xz)$, showing that $\deg f$ does not necessarily equal $\deg V(f)$. 
\end{ex}

\begin{ex}
If $f$ is any monomial, then $\Delta$ consists of a single point. Hence $V(f)$ has degree 0. This is appropriate since $V(f)$ is an empty set (Remark \ref{monom}).
\end{ex}
A vertex $V$ of a tropical curve is called {\em 3-valent} if $V$ has exactly 3 adjacent edges. Furthermore, if these edges have weights $m_1, m_2, m_3$ and primitive integer direction vectors $\feit{u}=(u_0,u_1,u_2), \feit{v}=(v_0,v_1,v_2),\feit{w}=(w_0,w_1,w_2)$ respectively, we define the {\em multiplicity} of $V$ to be the absolute value of the number
\begin{equation*}
  m_1m_2\left| \begin{array}{ccc}
u_0 & u_1 & u_2 \\
v_0 & v_1 & v_2 \\
1 & 1 & 1 \end{array} \right|=
m_2m_3\left| \begin{array}{ccc}
v_0 & v_1 & v_2 \\
w_0 & w_1 & w_2 \\
1 & 1 & 1 \end{array} \right|=
m_1m_3\left| \begin{array}{ccc}
u_0 & u_1 & u_2 \\
w_0 & w_1 & w_2 \\
1 & 1 & 1 \end{array} \right|.
\end{equation*}

\begin{remark}
This is a projective version of \cite[Definition 2.16]{Mikh}. Notice that if we consider the curve in the affine tropical plane $z=0$, the multiplicity of $V$ becomes $m_1m_2\cdot\text{Area}(\feit{u},\feit{v})$, which is in agreement with \cite[Definition 2.16]{Mikh}.
\end{remark}

\begin{Def}
A tropical curve is called {\em smooth} if every vertex is 3-valent and has multiplicity 1. 
\end{Def}

Notice that in a smooth tropical curve, every edge has weight 1.
We conclude this subsection by defining the {\em genus} of a smooth tropical curve:
\begin{Def}\label{genusdef}
Let $C=V(f)$ be a smooth tropical curve. The {\em genus} of $C$ is the number of interior lattice points of $\Subdiv_f$.
\end{Def}

\subsection{Intersections of tropical curves}
We say that two tropical curves $C$ and $D$ intersect {\em transversally} if no vertex of $C$ lies on $D$ and vice versa. In a transversal intersection we define intersection multiplicities as follows: Let $P$ be an intersection point of $C$ and $D$, where the two edges meeting have weights $m_1$ and $m_2$, and primitive integer direction vectors $(v_0,v_1,v_2)$ and $(w_0,w_1,w_2)$ respectively. The {\em intersection multiplicity} $\mu_P$ of $C$ and $D$ at $P$ is then the absolute value of
\begin{equation*}
  m_1m_2
\left| \begin{array}{ccc}
v_0 & v_1 & v_2 \\
w_0 & w_1 & w_2 \\
1 & 1 & 1 \end{array} \right|.
\end{equation*}

Non-transversal intersections are dealt with in the following way: For {\em any} intersecting tropical curves $C$ and $D$, let $C_\epsilon$ and $D_\epsilon$ be nearby translations of $C$ and $D$ such that $C_\epsilon$ and $D_\epsilon$ intersect transversally. We then have (\cite[Theorem 4.3]{RGST}):
\begin{theodef}\label{stable}
Let the {\em stable intersection} of $C$ and $D$, denoted $C\cap_{st} D$, be defined by
\begin{equation*}
C\cap_{st} D=\lim_{\epsilon\pil 0}(C_\epsilon \cap D_\epsilon).  
\end{equation*}
This limit is independent of the choice of perturbations, and is a well-defined subset of points with multiplicities in $C\cap D$.
\end{theodef}

\begin{theo}[Tropical Bezout]\label{Bez}
Assume $C$ and $D$ are tropical curves of degrees $c$ and $d$ respectively. If both curves have full support, then 
their stable intersection consists of $cd$ points, counting multiplicities.
\end{theo}

\begin{proof}
See \cite[Theorem 4.2 and Corollary 4.4]{RGST}. The idea is to show that the number of (stable) intersection points is invariant under translations of the curves. Thus we can arrange the two curves such that for each of them, the intersection points lie on the unbounded rays in one of the three coordinate directions. It is then trivial to check that $\sharp (C\cap_{st} D)=cd$. 
\end{proof}

\begin{remark}\label{bezrem}
In \cite{RGST}, the tropical semiring is defined as $(\rr\cup\{\infty\},\min,+)$ instead of $(\rr,\max,+)$ as here. The inclusion of the additive identity element $\infty$ makes $\pp_{tr}^2$ strictly larger than $\rr^2$, opening for the possibility that two tropical curves could have intersection points at infinity (like classical algebraic curves). However, it is not hard to see that two tropical curves, of which at least one has full support, always have all their intersection points in $\rr^2$. Hence the theorem holds with our definition as well.
\end{remark}

There is also a tropical version of Bernstein's theorem: Recall that the {\em mixed area} of two convex polygons $R$ and $S$ is defined as the number $\Area(R+S)-\Area(R)-\Area(S)$, where $R+S$ is the {\em Minkowski sum} of $R$ and $S$.

\begin{theo}[Tropical Bernstein]\label{Bern}
Let $C=V(f)$ and $D=V(g)$ be any tropical curves intersecting transversally, with Newton polygons $\Delta_f$ and $\Delta_g$ respectively. Then the number of intersection points, counting multiplicities, equals the mixed area of $\Delta_f$ and $\Delta_g$.
\end{theo}
\begin{proof}
See \cite[Theorem 9.5]{Stur}.
\end{proof}

Although arguably not as enlightening as the homotopy argument given in \cite{RGST}, one can prove Theorem \ref{Bez} as a special case of Theorem \ref{Bern}. In fact, we can get a stronger result:

\begin{theo}[Strong version of Tropical Bezout]\label{MBez}
Assume $C$ and $D$ are tropical curves of degrees $c$ and $d$ respectively. If at least one of the curves have full support, then their stable intersection consists of $cd$ points, counting multiplicities.
\end{theo}
\begin{proof}
Because of Theorem-Definition \ref{stable} we can assume the intersection is transversal. Note that for any positive integers $c$ and $d$, we have the Minkowski sum $\Ga_c+\Ga_d=\Ga_{c+d}$. Hence the mixed area of $\Ga_c+\Ga_d$ equals $\frac12(c+d)^2-\frac12 c^2-\frac12 d^2=cd$. This proves Theorem \ref{Bez}.

Suppose now $C$ has full support, i.e. $\Delta_f=\Ga_c$, and that $\Delta_g$ is a convex polygon
of the form $\Ga_d\smallsetminus Q$, where $Q\sub\Ga_d$ is a lattice polygon
containing exactly one of the corners of $\Ga_d$, say $(d,0,0)$. Then $\Area(\Delta_f+\Delta_g)=\Area(\Ga_c+(\Ga_d\smallsetminus Q))=\Area(\Ga_c+\Ga_d)-\Area(Q)$. Thus the mixed area of $\Delta_f$ and $\Delta_g$ is
\begin{multline*}
\Area(\Delta_f+\Delta_g)-\Area(\Delta_f)-\Area(\Delta_g)=\\
(\Area(\Ga_c+\Ga_d))-\Area(Q)-\Area(\Ga_d)-(\Area(\Ga_d)-\Area(Q))=cd.
\end{multline*} 

The same argument shows that we can do the same at the other corners, without changing the mixed area. In this way we can form
any Newton polygon
$\Delta_g$ associated to a tropical curve of degree $d$. Hence $\sharp(C\cap_{st} D)=cd$ for any tropical curve $D$ of degree $d$.
\end{proof}
 
\begin{ex}
If neither of the two curves have full support, the theorem will not hold in general. For example, if $C$ and $D$ are the quadric curves given by $C=V(x^2\oplus y)$ and $D=V(x\oplus y^2)$, then $C\cap D$ consists of a single point with multiplicity 3. Another example is given by the non-intersecting lines $V(0\oplus x)$ and $V(1\oplus x)$.
\end{ex}

An important special case of Theorem \ref{MBez} is the following corollary: 
\begin{cor}
Let $D$ be any tropical curve of degree $d$. Then any tropical line meets $D$ stably in exactly $d$ points, counting multiplicities.
\end{cor} 

\section{Divisors on smooth tropical curves}

Let $C$ be a smooth tropical curve in $\pp_{tr}^2$. 

\begin{Def}
We define the {\em group of divisors} on $C$, $\Div(C)$, to be the free abelian group generated by the points on $C$. A {\em divisor} $D$ on $C$ is an element of $\Div(C)$, i.e. a finite formal sum of the form
$D=\sum \mu_P P$.
\end{Def} 
 The number $\sum \mu_P$ is as usual called the {\em degree} of $D$. Observe that the elements of degree 0 in $\Div(C)$ form a group, denoted by $\Div^0(C)$.

\begin{Def}
Given a homogeneous tropical polynomial $f$, we define the associated divisor $\divv f\in\Div(C)$ as the formal sum of points in $C\cap_{st} V(f)$, counted with their respective intersection multiplicities.
The {\em principal divisors} on $C$ are the divisors of the form
\begin{equation*}
  \divv\frac{f}{g}:=\divv f-\divv g,
\end{equation*}
where $f$ and $g$ are homogeneous tropical polynomials of the same degree.
\end{Def}
\begin{Def}
Two divisors $D_1$ and $D_2$ are {\em linearly equivalent}, denoted as $D_1\sim D_2$, if $D_1-D_2$ is principal.  
\end{Def}

Linear equivalence is an equivalence relation, and as in the classical case one can show that it restricts to an equivalence relation on the subgroup $\Div^0(C)$. Hence we can make the following definition:
\begin{Def}
The group $\Div^0(C)/\nsp\sim$ is called the {\em Jacobian} of $C$, $\Jac(C)$.
\end{Def}

\begin{remark}
Most of the definitions in this section make sense also for non-smooth tropical curves. Imitating the Cartier divisors in classical algebraic geometry, it would be natural then for $\Div(C)$ to be generated by the {\em weighted points} on $C$, i.e. the set $\{m_P P\}$ where $m_P$ is the weight of the edge containing $P$. The main problem with this arises if $P$ is a vertex. As an example, we pose the following problem: What should be the ''weight'' of $P$ if $P$ is the vertex of the curve $V(z^8\oplus x^2y^6\oplus x^5z^3)$?
\end{remark}

\section{The group law on a tropical elliptic curve}
\subsection{Tropical elliptic curves in $\rr^2$.}

\begin{Def}
A {\em tropical elliptic curve} is a smooth tropical curve of degree 3 and genus 1.
\end{Def}

To make the notation simpler (and illustrations more effective), we will from now on regard our tropical curves as subsets of $\rr^2$ instead of $\pp^2_{tr}$. More precisely, we choose the affine plane $\rr^2\sub\pp_{tr}^2$ given by $z=0$. As commented in Remark \ref{ekstra} we don't lose any points on $C$ by doing this. Moreover, the strong version of Bezout's theorem, as stated in Theorem \ref{MBez}, holds for tropical curves in $\rr^2$ (cf. Remark \ref{bezrem}).

Thus in the following we assume that $C=V(f)\sub\rr^2$ is an affine tropical elliptic curve, where $f(x,y)=F(x,y,0)$ is the tropical dehomogenization of a homogeneous cubic tropical polynomial $F(x,y,z)$. The Newton polygon $\Delta_f$ is then contained in the triangle in $\zz^2$ with vertices $(0,0),(3,0),(0,3)$, and the associated subdivision $\Subdiv_f$ is a triangulation of $\Delta_f$ (since $C$ is smooth). The condition that $C$ is elliptic is by Definition \ref{genusdef} equivalent to $(1,1)$ being an interior point of $\Subdiv_f$. Hence $C$ contains a unique cycle, which we will denote by $\wb{C}$. Furthermore, we call each connected component of $C\smallsetminus \wb{C}$ a {\em tentacle} of $C$.

\begin{prop}
Let $P$ and $Q$ be points on the same tentacle of $C$. Then $P\sim Q$.
\end{prop}
\begin{figure}[htbp]
\begin{center}
\input{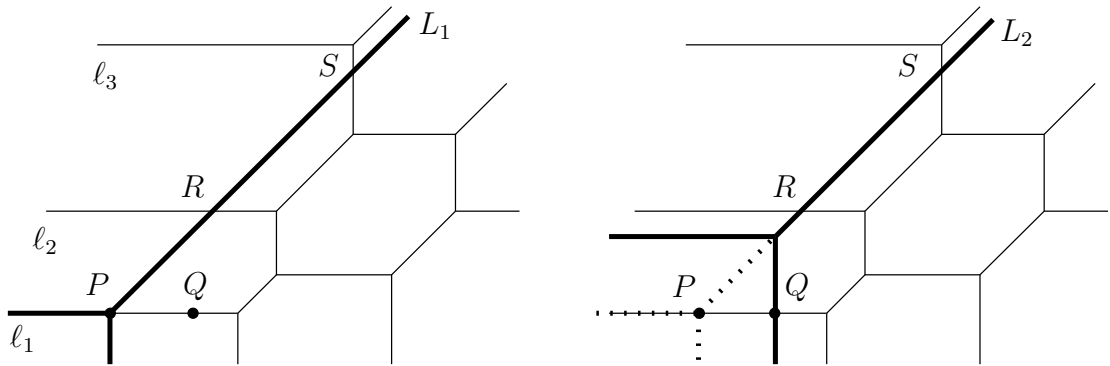}
\caption{Showing that points on $\ell_1$ are linearly equivalent.}\label{tent}
\end{center}
\end{figure}
\begin{proof}
We begin by showing that the points on any unbounded ray are equivalent. By symmetry, it is enough to prove this for the rays that are unbounded in, say, the $x$-coordinate. Figure \ref{tent} shows a typical situation with three such rays, $\ell_1$, $\ell_2$ and $\ell_3$.

The following argument shows that any two sufficiently close points $P$ and $Q$ on $\ell_1$ are equivalent: Assume $P$ is further away from $\wb{C}$ than $Q$, and consider a linear form $f_1$ such that $L_1=V(f_1)$ is the tropical line with center in $P$. Then $\divv f_1=P+R+S$, where $R$ and $S$ lies on the $z$-ray of $L_1$ (i.e. the ray with direction vector $(1,1)$). Now let $f_2$ be such that $L_2=V(f_2)$ is the line passing through $Q$ and with center on the $z$-ray of $L_1$. Then $\divv f_2=Q+R+S$ (as long as $P$ and $Q$ are close enough). It follows that $P-Q=\divv \frac{f_1}{f_2}$, in other words $P\sim Q$. 

To show that {\em any} two points $P$ and $Q$ on $\ell_1$ are equivalent, we can choose a finite sequence of points $P=P_1,P_2,\dotsc,P_m=Q$ on $\ell_1$ such that each pair $(P_i,P_{i+1})$ is close enough for the above technique to work. Then $P=P_1\sim\dotsb\sim P_m=Q$. 

\begin{figure}[htbp]
\begin{center}
\input{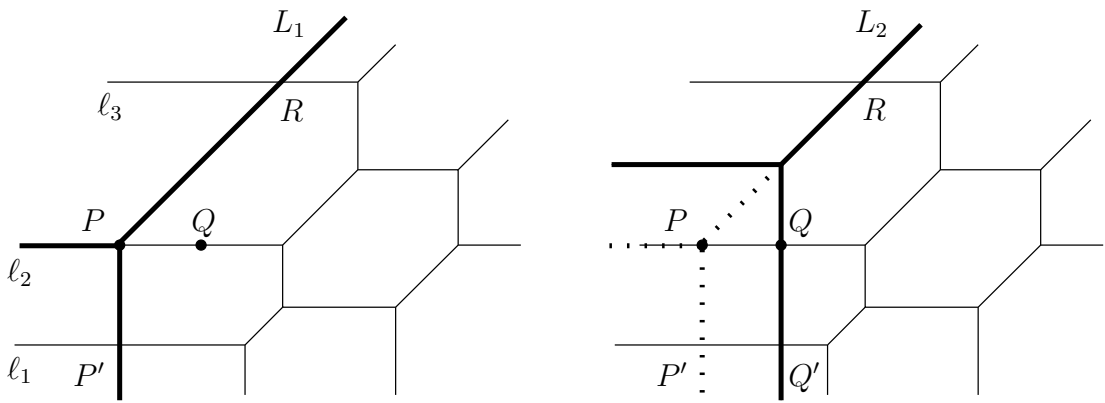}
\caption{Showing that $P\sim Q$ on $\ell_2$.}\label{tent2}
\end{center}
\end{figure}

A similar argument shows that the points on $\ell_2$ are equivalent. The idea is sketched in figure \ref{tent2}. To show that $P$ and $Q$ are equivalent, take the tropical line $L_1$ with center in $P$ and slide it along the $z$-ray (i.e. keeping $R$ as intersection point with $C$) until it passes through $Q$. With the notation on figure \ref{tent2}, we see that $P+P'+R\sim Q+Q'+R$. But $P'\sim Q'$, since they are on $\ell_1$, thus $P\sim Q$.

The same technique works for $\ell_3$ and also for the bounded line segments on the tentacles. It is easy to see that any tentacle of a tropical elliptic curve can be handled in this way. 
\end{proof}

\subsection{An explicit homeomorphism $\wb{C}\pil S^1$}
As before, let $C\sub \rr^2$ be a tropical elliptic curve, and let $\oo\in\wb{C}$ be a fixed point. Let $V_1,\dotsc,V_n$ be the vertices of $\wb{C}$ in counter-clockwise direction, such that if $\oo$ is a corner then $V_1=\oo$, otherwise $\oo$ lies between $V_1$ and $V_n$. Let $E_1,\dotsc, E_n$ be the edges of $\wb{C}$, such that $E_1=[V_1V_2]$ and so on. Let $|\;|$ denote Euclidean length.

Obviously, as a topological space, $\wb{C}$ is homeomorphic to $S^1$. We will now construct one such homeomorphism, piecewise-linear in the Euclidean metric.

For $i=1,\dotsc,n$, let $\e_i:=\frac{1}{|\feit{v}_i|}$, where $\feit{v}_i$ is the primitive integer vector along $E_i$. We define the {\em lattice length} of $E_i$ to be $\e_i|E_i|$. Observe that if $E_i$ has integral endpoints, then
$$\e_i|E_i|=1+\sharp\{\text{interior lattice points on $E_i$}\}.$$
Let $L$ be the total lattice length of $\wb{C}$, i.e.
$L=\e_1|E_1|+\dotsb+\e_n|E_n|$.

\begin{remark}
This gives $\wb{C}$ the so-called {\em $\zz$-affine structure} described by Mikhalkin in \cite[Remark 2.4]{Mikh}.
\end{remark}

We now define a homeomorphism $\la\colon \wb{C}\longrightarrow \rr/\zz\approx S^1$, linear in the Euclidean metric of each edge $E_i$. It is enough to specify the images of $\oo$ and the vertices, which we do recursively:
\begin{equation}\label{ladef}
  \begin{split}
    \la(\oo) & = 0\\
\la(V_1) & = \frac{\e_n\cdot |\oo V_1|}{L}\\
\la(V_{i+1}) & = \la(V_{i})+\frac{\e_{i}\cdot |E_i|}{L},\quad i=1,\dotsc,n-1.
  \end{split}
\end{equation}
Identifying $\rr/\zz$ with the interval $[0,1)$, we define the {\em signed lattice distance} $d_C(P,Q)$ between points $P$ and $Q$ on $\wb{C}$ by the formula
\begin{equation}\label{d-la}
  d_C(P,Q)=L\cdot(\la(Q)-\la(P)).
\end{equation}
Notice that for any three points $P,Q,R\in\wb{C}$ we have
\begin{equation*}
d_C(P,Q)+d_C(Q,R)=d_C(P,R).  
\end{equation*}

\subsection{Description of the group law}
We now move on to the task of determining when divisors of the form $P+Q$ are linearly equivalent. When trying to imitate the techniques from the classical case, we stumble across the following problem:
Given two points $P$ and $Q$ on $\wb{C}$, we cannot always find a tropical line $L$ that intersects $C$ stably in $P$ and $Q$. (Recall that a stable intersection is defined as a limit of transversal intersections.) If there exists such a tropical line, we call $(P,Q)$ a {\em good pair}. 

We fix notation $\feit{p}_1=(-1,0)$, $\feit{p}_2=(0,-1)$ and $\feit{p}_3=(1,1)$ for the primitive integer direction vectors of a tropical line.

\begin{lem}\label{sumlem}
Let $P,Q,P',Q'$ be any points on $\wb{C}$. Then 
$$P+Q\sim P'+Q'\quad\Longleftrightarrow\quad d_C(P,P')=-d_C(Q,Q').$$ 
\end{lem}

\begin{proof}
We proceed in two steps. First, we prove the result when $(P,Q)$ and $(P',Q')$ are good pairs. Using this, we then generalize to any pairs.

$\bullet$ {\em Step 1}.
Assume $(P,Q)$ and $(P',Q')$ are good pairs, and that $P+Q\sim P'+Q'$. Then there exists (unique) tropical lines $L$ and $L'$, and a point $R\in\wb{C}$ such that $L\cap_{st} C=P+Q+R$ and $L'\cap_{st} C=P'+Q'+R$. Consider a homotopy $L_t$ of lines containing $R$ such that $L_0=L$ and $L_1=L'$. It is enough to consider the case where $P$ and $P'$, and $Q$ and $Q'$, are on the same edge respectively, and where $L'$ is a parallell displacement of $L$ along one of the axes. Indeed, in more complex cases, the homotopy can be broken down into parts with the above properties. 

\begin{figure}[tbp]
\noindent
\begin{minipage}[b]{.48\linewidth}
\begin{center}
\input{move.pstex_t}
\caption{Illustrating Step 1.}\label{move}
\end{center}
\end{minipage}\hfill%
\begin{minipage}[b]{.48\linewidth}
\begin{center}
\input{proj.pstex_t}
\caption{Non-orthogonal projection}\label{proj}
\end{center}
\end{minipage}\end{figure}

Let $\feit{v}_P$ and $\feit{v}_Q$ be the primitive integer vectors corresponding to the edges of $\wb{C}$ containing $P$ and $Q$ (see figure \ref{move}). Now assume that $L'$ equals the shifting of $L$ $\delta$ units in the direction of, say, $\feit{p}_1$. Then from the general formula for (non-orthogonal) vector projection (figure \ref{proj}), we find the displacements of $P$ and $Q$: 

\begin{equation}\label{euc}
\begin{split}
 PP'&= \,\frac{|\feit{p}_2\times \delta\feit{p}_1|}{|\feit{p}_2\times \feit{v}_P|}\:\feit{v}_P\,=
\delta\feit{v}_P\quad\Longrightarrow\quad |d_C(P,P')|=\frac1{|\feit{v}_P|}|\delta\feit{v}_P|=\delta,\\
QQ'&= \,\frac{|\feit{p}_3\times \delta\feit{p}_1|}{|\feit{p}_3\times \feit{v}_Q|}\:\feit{v}_Q\,=
\delta\feit{v}_Q\quad\Longrightarrow\quad |d_C(Q,Q')|=\frac1{|\feit{v}_Q|}|\delta\feit{v}_Q|=\delta.
\end{split}
\end{equation}
(Notice that both the denominators above equals 1, since the intersections at hand have multiplicity 1.)
According to the orientation of $\wb{C}$, $P$ and $Q$ are moved in opposite direction. Hence $d_C(P,P')=-d_C(Q,Q')$ as claimed.

The implication $\Leftarrow$ follows by a similar argument.

$\bullet$ {\em Step 2}.
Now assume $(P,Q)$ is not a good pair. Let $L_1$ and $L_2$ be tropical lines through $P$ and $Q$ respectively, and let $R_1,S_1,R_2,S_2$ be the other intersection points. The idea is to move $L_1$ and $L_2$ into new lines $L_1'$ and $L_2'$ in such a way that $R_1,S_1,R_2,S_2$ are preserved as intersection points. $P$ and $Q$ will not be preserved; they will move to new points $P'$ and $Q'$. (See figure \ref{move2}.) By construction, these points satisfy $P'+Q'\sim P+Q$. Using our results in Step 1 on each of the lines $L_1$ and $L_2$, it follows that $d_C(P,P')=-d_C(Q,Q')$. Conversely, it is not hard to see that in this way one can reach any nearby pair $(P',Q')$ satisfying $d_C(P,P')=-d_C(Q,Q')$. 

Finally, by chosing $L_1$ and $L_2$ wisely, $(P',Q')$ will form a good pair. Since we proved in Step 1 that the lemma is true for good pairs, it then follows that the lemma holds for any pairs $(P,Q)$ and $(P',Q')$.
\begin{figure}[htbp]
\begin{center}
\input{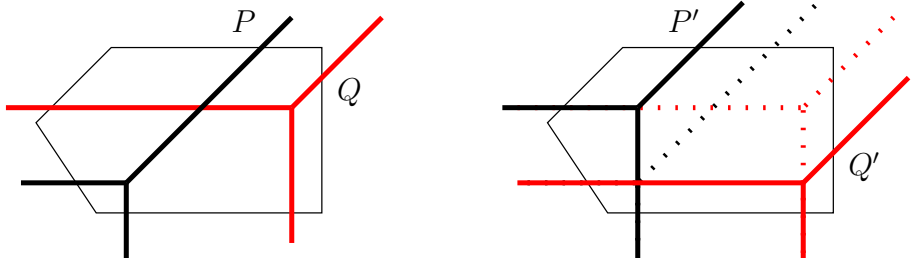}
\caption{Moving a bad pair $(P,Q)$.}\label{move2}
\end{center}
\end{figure}
\end{proof}

\begin{prop}\label{bij}
For any fixed point $\oo\in\wb{C}$, 
the map $\tau_\oo\colon \wb{C}\pil \Jac(C)$ given by $P\mapsto P-\oo$ is a bijection of sets.
\end{prop}

\begin{proof}
Injectivity follows immediately from Lemma \ref{sumlem}, since
\begin{equation*}
P-\oo\sim Q-\oo\Longrightarrow  P+\oo\sim Q+\oo \Longrightarrow d_C(P,Q)=0\Longrightarrow P=Q.
\end{equation*}

To prove surjectivity, let $D$ be any divisor of degree 0 . We must show that there exists $P\in \wb{C}$ such that $D\sim P-\oo$. Assume first that $D=P_1-Q_1$, where $P_1,Q_1\in \wb{C}$. Choose $P$ such that $d_C(P,P_1)=d_C(\oo,Q_1)$, then Lemma \ref{sumlem} gives $P+Q_1\sim P_1+\oo$. Thus we have $D=P_1-Q_1\sim P-\oo$.

Now assume $D=D_1-D_2$, where $D_1=P_1+\dotsb+P_n$ and $D_2=Q_1\dotsb+Q_n$ are any effective divisors of degree $n>1$. Let $P_{12}$ and $Q_{12}$ be points such that $P_1+P_2\sim \oo+P_{12}$ and $Q_1+Q_2\sim \oo+Q_{12}$. Then
\begin{equation*}
D\sim \oo +P_{12}+\dotsb+P_n-(\oo+Q_{12}+\dotsb+Q_n)=P_{12}+\dotsb+P_n-(Q_{12}+\dotsb+Q_n).
\end{equation*}
Hence $D\sim D_1'-D_2'$, where $D_1$ and $D_2$ are effective of degree $n-1$. This way we can reduce to the case $n=1$, which we already proved. 
\end{proof}

Because of Proposition \ref{bij}, $\wb{C}$ has a natural group structure:
\begin{Def}
Define $(\wb{C},\oo)$ to be the group consisting of points on $\wb{C}$, with the group structure induced from $\Jac(C)$ such that the bijection $\tau_\oo$ is an isomorphism of groups.
\end{Def}

The next theorem and its corollary are the main results of this paper.
\begin{theo}
Let $P$ and $Q$ be any points on $\wb{C}$, and let $\feit{+}$ denote addition in the group $(\wb{C},\oo)$. Then the point $P\feit{+}Q$ satisfies the relation $$d_C(\oo,P\feit{+}Q)=d_C(\oo,P)+d_C(\oo,Q).$$
\end{theo}
\begin{proof}
To simplify notation, let $R=P\feit{+}Q$. Then because $\tau_\oo$ is a group isomorphism, we have $$R-\oo=\tau_\oo(R)=\tau_\oo(P)+\tau_\oo(Q)=P-\oo+Q-\oo.$$ Thus $R\sim P+Q-\oo$, i.e. $R+\oo\sim P+Q$, and it follows from Lemma \ref{sumlem} that
\begin{equation*}
  d_C(P,R)=d_C(\oo,Q).
\end{equation*}
Adding $d_C(\oo,P)$ on each side then gives $d_C(\oo,R)=d_C(\oo,P)+d_C(\oo,Q)$ as wanted.
\end{proof}

\begin{remark}
We can describe the group law geometrically just as in the classical case of elliptic curves: To add $P$ and $Q$ we do the following. If $(P,Q)$ is a good pair, consider the tropical line $L$ through $P$ and $Q$, and let $R$ be the third intersection point of $L$ and $\wb{C}$. Now if $(R,\oo)$ is a good pair, let $L'$ be the through $R$ and $\oo$. Then $P\feit{+}Q$ is the third intersection point of $L'$ and $\wb{C}$. (See figure \ref{add} for an example.)

\begin{figure}[htbp]
\begin{center}
\input{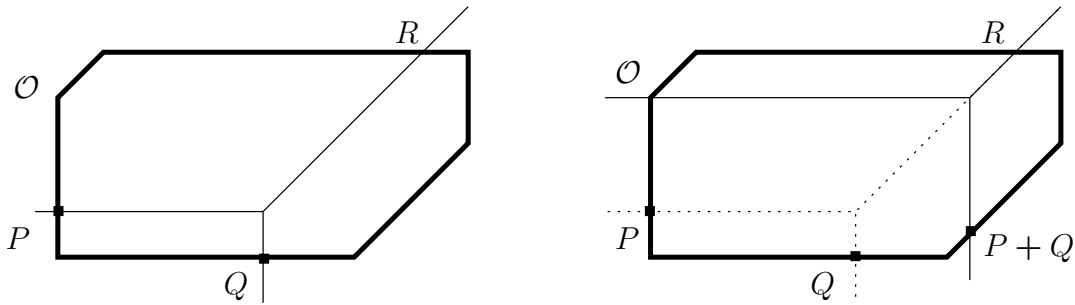}
\caption{Adding points on a tropical elliptic curve.}\label{add}
\end{center}
\end{figure}

If any of the pairs $(P,Q)$ and $(R,\oo)$ fail to be good, then move the two points involved equally far (in the lattice metric) in opposite directions until they form a good pair, and use this new pair as described above.
\end{remark}

\begin{cor}
The map $\la\colon (\wb{C},\oo)\longrightarrow \rr/\zz\approx S^1$ defined in \eqref{ladef} is a group isomorphism.
\end{cor}

\begin{proof}
It follows from the relation \eqref{d-la} that for any $P$ we have $\la(P)=d_C(\oo,P)/L$. Thus
  \begin{equation*}
\la(P\feit{+}Q)=\frac{d_C(\oo,P\feit{+}Q)}{L}=\frac{d_C(\oo,P)}{L}+\frac{d_C(\oo,Q)}{L}=\la(P)+\la(Q).
  \end{equation*}
\end{proof}

\noindent
{\it Acknowledgements}.
Thanks to Grigory Mikhalkin for helpful advice at an early stage of this work, and to Kristian Ranestad for many useful discussions, suggestions and corrections.

\bibliographystyle{plain}
\bibliography{trop_ref}
\end{document}